\documentclass[11pt]{article}

\usepackage{amsmath} 
\usepackage{amssymb}
\usepackage{amsthm}
\usepackage{mathtools}
\usepackage{fancyhdr}
\usepackage{hyperref}
\usepackage{cleveref}
\usepackage[square, numbers]{natbib}
\usepackage[english]{babel}
\usepackage[utf8]{inputenc}
\usepackage{tikz}
\usepackage{enumitem}

\newtheorem{theorem}{Theorem}[section]
\newtheorem*{theorem-nonum}{Theorem}
\newtheorem{definition}[theorem]{Definition}
\newtheorem*{definition-nonum}{Definition}
\newtheorem{lemma}[theorem]{Lemma}

\newtheorem{remark}[theorem]{Remark}
\newtheorem{corollary}[theorem]{Corollary}
\newtheorem{example}[theorem]{Example}
\newtheorem{proposition}[theorem]{Proposition}

\title{Characterization of minimizers of an anisotropic variant of the Rudin-Osher-Fatemi functional with $L^1$ fidelity term}
\author{Nils Dabrock\thanks{e-mail: nils.dabrock@tu-dortmund.de}\\
Fakultät für Mathematik \\ Technische Universität Dortmund}

\newcommand{\R}{{\mathbb R}}%
\newcommand{\N}{{\mathbb N}}%

\newcommand\restr[2]{{
  \left.\kern-\nulldelimiterspace 
  #1 
  \vphantom{\big|} 
  \right|_{#2} 
  }}

\newcommand{\divergence}{\nabla \cdot}

\DeclareMathOperator{\argmax}{argmax}

\DeclareMathOperator{\TV}{TV}
\DeclareMathOperator{\support}{spt}
\DeclareMathOperator{\interior}{int}

\newcommand{\Mn}{[\mathcal{M}(\R^n)]^n}
\newcommand{\symmdiff}{\,\triangle \,}
\newcommand{\rmid}{\,\middle| \,}

\begin{document}
	\bibliographystyle{siam}
	\maketitle

	\begin{abstract}
		In this paper we study an anisotropic variant of the Rudin-Osher-Fatemi functional with $L^1$ fidelity term of the form
		\[
			E(u) = \int_{\R^n} \phi(\nabla u) + \lambda \| u -f \|_{L^1(\R^n)}.
		\]
		We will characterize the minimizers of $E$ in terms of the Wulff shape of $\phi$ and the dual anisotropy.
		In particular we will calculate the subdifferential of $E$.
		We will apply this characterization to the special case $\phi = |\cdot|_1$ and $n=2$, which has been used in the denoising of 2D bar codes.
		In this case, we determine the shape of a minimizer $u$ when $f$ is the characteristic function of a circle.
	\end{abstract}	
	
	\section{Introduction}
		We analyze an anisotropic variant of the Rudin-Osher-Fatemi functional (anisotropic $\TV$-$L^1$ energy in the following)
		\[
			E(u) = \int_{\R^n} \phi(\nabla u) + \lambda \| u -f \|_{L^1(\R^n)},
		\]
		which has been proposed by Choksi et al.~in \cite{Choksi2011} for the special case of $\phi = |\cdot|_1$ and $n=2$ for the denoising of 2D bar codes. Given a noisy input $f \in L^1(\R^n)$ we seek
		for a minimizer $u$ of $E$. 
		
		The famous Rudin-Osher-Fatemi functional \cite{Rudin1992} is given by
		\[
			E_\text{ROF}(u) = \int_{\Omega} |\nabla u|_2 + \lambda \|u-f\|_{L^2(\Omega)}^2,
		\]
		where the noisy input is modeled by a function $f \in L^2(\Omega)$ typically defined on a bounded domain $\Omega \subset \R^2$. 
		The first term measures the regularity of $u$ with respect to the Euclidean metric on $\R^2$. We call this term \textit{isotropic total variation}.
		The second term measures the distance to the original signal and is called \textit{fidelity term}.
		This model is not contrast invariant: plugging in a rescaled noisy image $cf$ for some rescaling parameter $c > 0$ does not result in a rescaled minimizer $cu$. Furthermore the minimization of $E_\text{ROF}$ is not faithful to any $f \ne 0$, which means that there is no $f \ne 0$ such that $f$ is the minimizer of $E_\text{ROF}$. 
		
		In \cite{Chan2005} Chan and Esedo\=glu analyze the isotropic $\TV$-$L^1$ model. This model has a different fidelity term, where the squared $L^2$ distance to the measured signal $f$ is replaced by the $L^1$ distance to $f$. The authors prove that this modification yields a contrast invariant minimization problem. They show that $u$ is a minimizer if and only if the level sets of $u$ solve a shape optimization problem. This motivates the study of the geometric properties of minimizers. They can show that the $\TV$-$L^1$ model is faithful to every characteristic function of a bounded domain with $C^2$ boundary, if $\lambda$ is larger than some domain dependent constant. A negative result is, that the characteristic function of a square can not appear as a minimizer of this model. 
		
		Therefore, when working with images that consist of shapes with some anisotropic feature, a modified total variation must be considered. 
		This has been done by Esedo\=glu and Osher in \cite{Esedoglu2004}
		for the classical Rudin-Osher-Fatemi model. They consider an energy 
		\[
			E_\text{EO}(u) = \int_{\Omega} \phi(\nabla u) + \lambda \|u-f\|_{L^2(\Omega)}^2
		\]
		where $\phi$ is a norm on $\R^2$. They prove that this energy prefers the corresponding Wulff shape
		\[
			W_\phi = \left\{ x \in \R^2 \rmid x\cdot y \le \phi(y) ~\forall y \in \R^2 \right\}.
		\]
		In particular they show that for large $\lambda$ the unique minimizer for $f = \chi_{W_\phi}$ is given by $cf$ for some constant $c > 0$.
		
		In this work we combine the techniques used in \cite{Chan2005} and \cite{Esedoglu2004} to analyze the functional $E$. In his PhD thesis \cite{Duval2011} Duval starts with an excellent overview over the $\TV$-$L^1$ model. He extends the geometric point of view of \cite{Chan2005} and proposes two algorithms based on these results. At the end of the first part of \cite{Duval2011} he mentions
		the generalization to anisotropies described by a norm $\phi$. If in addition $\phi$ is crystalline, that means that the Wulff shape is a polytope, he characterizes the minimizers for $f = \chi_C$ where $C$ is a bounded convex set. In this work we will allow a slightly more general class of anisotropies described by a Finsler metric $\phi$.
		
		The main result in this work is \Cref{thm:main}, which gives a dual characterization of minimizers $u \in BV(\R^n)$. We will use this result in \Cref{ex} to show that for $\phi = |\cdot|_1$
		and $f = \chi_{\{ |x|_2 \le 1 \}}$ the minimizer is given as intersection of the circle with a properly scaled square.
		In \cite{Choksi2011} Choksi et al. prove one implication of \Cref{thm:main} for this special case and deduce that $E$ is faithful to every 2D bar code as long as $\lambda$ is larger than some constant which depends on the size of the bar code. We conclude by \Cref{ex} that the converse conclusion is wrong, namely that there are faithful binary signals which are not given by 2D bar codes.
	
	\section{Anisotropic total variation}
		We consider an anisotropy given by a convex and nonnegative, positively 1-homogeneous function $\phi : \R^n \rightarrow [0, \infty)$,
		which means that $\phi$ satisfies $\phi(\alpha y) = \alpha \phi(y)$ for all $\alpha > 0$, $y \in \R^n$ and $\phi(y) = 0$ if and only if $y = 0$. A function $\phi$ like this is sometimes also called a gauge function  (see \cite{Freund1987}, \cite{Friedlander2014}) or a Finsler metric (see \cite{Kawohl2008}). The Wulff shape associated with $\phi$ is given by
		\[ W_\phi = \left\{ x \in \R^n  \rmid -x\cdot y \le \phi(y) ~\forall y \in \R^n \right\} \]
		and coincides with the polar set (sometimes also called the one-sided polar set) of
		\[ B_\phi = \left\{ y \in \R^n \rmid \phi(-y) \le 1\right\}. \]
		Then $W_\phi$ is a convex and compact set with $0 \in \interior W_\phi$.
		
		\begin{remark} Note that we need to introduce a minus sign in the definition of the Wulff shape, which differs from the classical setting. This is due to the fact that we want to allow anisotropies which are not even functions. This modification ensures that the Wulff shape is still the optimal shape for the anisotropic total variation we will define in \Cref{defn:aniso_tv}. This technical aspect occurs several times in this work. If $\phi$ is an even function, we can skip the minus sign.
		\end{remark}
		
		The gauge dual $\phi^\circ$ of $\phi$ is given by
		\[ \phi^\circ(x) \coloneqq \max_{\phi(y) \le 1} x\cdot y \]
		and can also be understood as the Minkowski function of $-W_\phi$,
		\[\phi^\circ(x) = \inf \{ \lambda > 0 \mid x \in -\lambda W_\phi \}. \]
		Since $\phi^{\circ\circ} = \phi$ (see \cite[Proposition 2.1]{Friedlander2014}), we have $B_\phi = W_{\phi^\circ}$ and $W_\phi = B_{\phi^\circ}$. The gauge functions $\phi$ and $\phi^\circ$ satisfy the Cauchy-Schwarz inequality in the sense that 
		\begin{equation}
			\label{eq:csi}
			 x\cdot y \le \phi^\circ(x) \phi(y) \text{ for all } x,y \in \R^n. 
		\end{equation}
		In the particular case $\phi = |\cdot|_p$ for some $p \in [1,\infty]$ we have $\phi^\circ = |\cdot|_q$, where $q$ is given by
		$\frac{1}{p} + \frac{1}{q} = 1$ in the usual sense. In this case the sets $W_\phi$ and $B_\phi$ are given by
		\begin{align*}
			W_\phi &= \{ x \in \R^n \mid |x|_q \le 1 \} \text{ and } \\
			B_\phi &= \{ x \in \R^n \mid |x|_p \le 1 \}.
		\end{align*}
		We need the following property, which is a generalization of the well known isometry of the embedding $L^\infty(\R^n ; \R^n) \hookrightarrow (L^1(\R^n ; \R^n))^*$ into the dual space of $L^1(\R^n ; \R^n)$.
		
		\begin{lemma}
			\label{eq:hoelder_extremal}
			Let $f \in L^\infty(\R^n ; \R^n)$ satisfy 
			\[ \int_{\R^n} f \cdot g \le \int_{\R^n} \phi(g) ~\forall g \in L^1(\R^n ; \R^n). \]
			Then $\phi^\circ(f) \le 1$ almost everywhere in $\R^n$, that means $f \in -W_\phi$ almost everywhere in $\R^n$.
			
			\begin{proof}
				Choose $\eta_{\phi^\circ}(x) \in \argmax_{\phi(y) \le 1} x \cdot y$ such that $\eta_{\phi^\circ}(x) \cdot x = \phi^\circ(x)$ and $\phi(\eta_{\phi^\circ}(x)) = 1$ for $x \in \R^n$. This can also be understood as $\eta_{\phi^\circ}(x) \in \partial \phi^\circ(x)$, where $\partial \phi^\circ(x)$ is the subdifferential of $\phi^\circ$ at $x \in \R^n$. 
				Let $g \in L^1(\R^n)$. Then, using \eqref{eq:csi}
				\begin{align*}
					\int_{\R^n} \phi^\circ(f) g &\le \int_{\R^n}  \phi^\circ(f) |g| \\
					&= \int_{\R^n} f \cdot \eta_{\phi^\circ}(f) |g|
					\le \int_{\R^n} \phi\left(\eta_{\phi^\circ}(f) |g|\right)
					= \int_{\R^n} |g|.
				\end{align*}
				Since $g$ was arbitrary we conclude $\phi^\circ(f) \le 1$ almost everywhere in $\R^n$.
			\end{proof}	
		\end{lemma}
		
		\begin{definition}[Anisotropic total variation, see \cite{Esedoglu2004}]
			\label{defn:aniso_tv}
			The anisotropic total variation of an $\R^n$-valued Radon measure $\mu \in \Mn$ is
			\[
			\int_{\R^n} \phi(\mathrm{d}\mu) \coloneqq \sup \left\{ \int_{\R^n} \varphi ~\mathrm{d}\mu \rmid \varphi \in C^1_c(\R^n ; \R^n), \varphi(x) \in -W_\phi ~ \forall x \in \R^n \right\}.
			\]
			The anisotropic total variation of a function $u \in BV(\R^n)$ is given by
			\[ \TV_\phi(u) \coloneqq \int_{\R^n} \phi(\mathrm{d} Du ) \]
			where $Du$ is the total variation measure of $u$. 
			
			\emph{Remark:} For $u \in BV(\R^n)$ we have 
			\[ \TV_\phi(u) = \sup \left\{ -\int_{\R^n} u \divergence \varphi \rmid \varphi \in C^1_c(\R^n ; \R^n), \varphi(x) \in -W_\phi ~ \forall x \in \R^n \right\}. \]
			The isotropic total variation is given by $\TV = \TV_{|\cdot|_2}$.
		\end{definition}

		The definition of the anisotropic total variation can be extended to functions that have a weak divergence in $L^\infty(\R^n)$.
		\begin{lemma}
			For $u \in BV(\R^n)$ we have
			\[
			\begin{aligned}
				\TV_\phi(u) = \sup \bigg\{ -\int_{\R^n} u \divergence v \,\bigg|\, &v \in L^\infty(\R^n ; \R^n), \\ 
				&\divergence v \in L^\infty(\R^n), v \in -W_\phi \text{ a.e.~in } \R^n \bigg\}. 
			\end{aligned}
			\]
			\begin{proof}
				The proof is an easy adaption of \cite[A.4]{Choksi2011}.
			\end{proof}
		\end{lemma}
		
		\begin{example}
			If $\Omega \subset \R^n$ is a bounded domain with Lipschitz continuous boundary, the anisotropic total variation equals
			\[ \TV_\phi(\Omega) \coloneqq \TV_\phi(\chi_\Omega) = \int_{\partial \Omega} \phi(-\nu). \]
			Here $\nu$ denotes the outward unit normal which is well defined almost everywhere on $\partial \Omega$.
			
			\begin{proof}
				We can prove this result by approximating the function 
				$\eta_\phi \circ (-\nu)$ (see \Cref{eq:hoelder_extremal}) by smooth functions and applying the Gauss formula.
			\end{proof}
		\end{example}
		
		\begin{corollary}
			\label{cor:tv_wulff}
			The total variation of the Wulff shape is given by
			\[
				\TV_\phi(W_\phi) = n |W_\phi|.
			\]
			\begin{proof}
				Since $W_\phi$ is convex, it has a Lipschitz continuous boundary, see \cite[1.2.2.3]{Grisvard1985}. Furthermore we know that at almost every point $x \in \partial W_\phi$ the Wulff shape lies on one side of the hyperplane $\{ y \in \R^n \mid y \cdot \nu(x) = x \cdot \nu(x) \}$. In particular we get 
				\[
					\phi(-\nu(x)) = \phi^{\circ\circ}(-\nu(x))
					= \sup_{\phi^\circ(y) \le 1} -\nu(x) \cdot y
					= \sup_{y \in B_{\phi^\circ}=W_\phi} \nu(x) \cdot y = x \cdot y. \]
				Then, using the Gauss formula
				\[ 
					\TV_\phi(W_\phi) = \int_{\partial W_\phi} \phi(-\nu) = \int_{\partial W_\phi} \nu(x) \cdot x ~\mathrm{d}x = n |W_\phi|.
				\]
			\end{proof}
		\end{corollary}
		
		\begin{proposition}
			\label{prop:equi}
			There are constants $c, C > 0$ satisfying
			\[ c  \TV_\phi(u) \le \TV(u) \le C \TV_{\phi}(u) ~\forall u \in BV(\R^n). \]
			
			\begin{proof}
				This follows directly from the fact that $W_\phi$ is compact and contains $0$ in its interior.
			\end{proof}
		\end{proposition}
		
		We state some results about the anisotropic total variation. All proofs are easy adaptions of the results \cite[5.2 and 5.5]{Evans1992} for the isotropic case.
		\begin{proposition}[Lower semi-continuity]
			\label{prop:lsc}
			Let $(u_k)_{k\in\N} \subset BV(\R^n)$, $u \in BV(\R^n)$ with $u_k \to u$ in $L^1_\text{loc}(\R^n)$, then $\TV_\phi(u) \le \liminf_{k\to\infty} \TV_\phi(u_k)$.
		\end{proposition}
		
		\begin{proposition}[Approximation]
			For every $u \in BV(\R^n)$ there is a sequence $(u_k)_{k \in \N} \subset BV(\R^n) \cap C^\infty(\R^n)$ satisfying 
			\begin{enumerate}
				\item $u_k \to u$ in $L^1(\R^n)$ and
				\item $\TV_\phi(u_k) \to \TV_\phi(u)$ for $k \to \infty$.
			\end{enumerate}
		\end{proposition}
		
		\begin{proposition}[Coarea formula]
			Let $u \in BV(\R^n)$. Define the level set $\Sigma(t) \coloneqq \{ x \in \R^n \mid u(x) > t \}$. Then we have
			\[ \TV_\phi(u) = \int_{-\infty}^\infty \TV_\phi(\Sigma(t)) ~\mathrm{d}t. \]
		\end{proposition}
		
	\section{Anisotropic $\TV$-$L^1$ model}
		\begin{definition}
			For a given function $f \in L^1(\R^n)$ and $\lambda > 0$ we consider the functional
			\[ E(u) \coloneqq E(u ; f, \lambda) \coloneqq \TV_\phi(u) + \lambda \| u - f \|_{L^1(\R^n)} ~,~ u \in BV(\R^n) \]
			and seek for a minimizer $u \in BV(\R^n)$ of $E$.
		\end{definition}
		
		Using \Cref{prop:equi}, the isotropic compactness result in \cite[5.2 Theorem 4]{Evans1992} and \Cref{prop:lsc}, we can apply the direct method of the Calculus of Variations to deduce the existence of a minimizer.
		
		Since $E$ is not strictly convex, the minimizer is not necessarily unique.
		
		\begin{example}
			Let $n=2$, $\phi = |\cdot|_2$, $B = B_\phi = \{ |x|_2 \le 1 \}$, $f = \chi_{B}$ and $\lambda = 2$. From what we will prove in this work, we infer that $f$ is a minimizer of $E$. For this special choice of $\lambda$ we note that  $\alpha f$ for $\alpha \in [0,1]$ is also a minimizer of $E$ because
			$E$ is convex and
			\[ E(f) = TV(f) = 2 \pi = 2\|f\|_{L^1} = E(0). \]
		\end{example}		
		
		For small $\lambda$ and functions $f$ with compact support, we only have the trivial minimizer.
		\begin{proposition}
			\label{prop:trivial_minimizer}
			Let $f \in L^1(\R^n)$ and $R > 0$ with $\support f \subset R\,W_\phi = \{ Rx \mid x \in W_\phi \}$. 
			The zero function is the unique minimizer of $E$ for all $0 < \lambda < \lambda_0 = \frac{n}{R}$.
			
			\begin{proof}
				This proof is similar to the proof of \cite[Lemma 2.2]{Choksi2011}, but we will additionally prove that $\lambda_0$ can be chosen as $\lambda_0 = \frac{n}{R}$.
				
				We know that the Wulff shape is the shape that minimizes the anisotropic total variation with prescribed area. 
				We refer to \cite{Fonseca1991a} for the proof. Using a scaling argument we can deduce that
				\begin{equation}
					\label{eq:iso_ineq}
					|A|^\frac{n-1}{n} \le C \TV_\phi(A)
				\end{equation}
				for all bounded sets $A \subset \R^n$ with finite perimeter. The constant $C$ is given by
				\[
					C = \frac{|W_\phi|^\frac{n-1}{n}}{\TV_\phi(W_\phi)} = n^{-1}|W_\phi|^{-\frac{1}{n}},
				\]
				where we applied \Cref{cor:tv_wulff}.
				Now let $u \in BV(\R^n)$ and $0 < \lambda < \lambda_0 = \frac{n}{R}$. We have that
				\begin{equation}
					\label{eq:splitting}
					\TV_\phi(u) = \TV_\phi(u_+) + \TV_\phi(-u_-)
				\end{equation}
				where $u_+$ and $u_-$ are the positive and negative parts of $u$. Since we are working with an anisotropy $\phi$ which is not necessarily an even function, we need to be careful
				when comparing the anisotropic total variation of positive and negative functions. We know that
				\[
					\int_{\R^n} -u \divergence \varphi = \int_{\R^n} u(x) \divergence (\varphi(-\cdot))(-x) ~\mathrm{d} x
					= \int_{\R^n} u(-x) \divergence (\varphi(-\cdot))(x) ~\mathrm{d}x
				\]
				for all $\varphi \in C^1_c(\R^n)$. Therefore the anisotropic total variation of the negative of $u$ is given by
				\begin{equation}
					\label{eq:negativetv}
					\TV_\phi(-u) = \TV_\phi(u(-\cdot)).
				\end{equation}
				Plugging \eqref{eq:negativetv} into \eqref{eq:splitting} and using the coarea formula together with the isoperimetric inequality \eqref{eq:iso_ineq} gives
				\begin{equation}
				\begin{aligned}
				\label{eq:tv_l1}
					\TV_\phi(u) &= \TV_\phi(u_+) + \TV_\phi(u_-(-\cdot)) \\
					&= \int_0^\infty \TV_\phi(\{u_+ > t\}) + TV_\phi(\{u_-(-\cdot) > t\}) ~\mathrm{d}t \\
					&\ge C^{-1} \int_0^\infty |\{u_+ > t\}|^\frac{n-1}{n} + |\{u_- > t\}|^\frac{n-1}{n} ~\mathrm{d}t \\
					&\ge C^{-1} \int_0^\infty |\{ |u| > t\}|^\frac{n-1}{n} ~\mathrm{d}t \\
					&\ge C^{-1} \int_0^\infty |\{ |u| > t\}\cap R\,W_\phi|^\frac{n-1}{n} ~\mathrm{d}t \\
					&\ge C^{-1} \int_0^\infty |\{ |u| > t\}\cap R\,W_\phi||R\,W_\phi|^{-\frac{1}{n}} ~\mathrm{d}t 
					= \lambda_0 \| u \|_{L^1(R\,W_\phi)}.
				\end{aligned}
				\end{equation}
				We conclude
				\begin{equation}
				\label{eq:e_zero}
				\begin{aligned}
					E(u) &\ge \lambda_0 \| u \|_{L^1(R\,W_\phi)} + \lambda \| u -f \|_{L^1(\R^n)} \\
					&\ge \lambda \|f\|_{L^1(\R^n)} = E(0).
				\end{aligned}
				\end{equation}
				It is easy to see, that the last inequality in \eqref{eq:e_zero} or the second to last inequality in \eqref{eq:tv_l1} is strict if $u \ne 0$.
			\end{proof}
		\end{proposition}
		
		\begin{remark}
			The constant $\lambda_0$ in \Cref{prop:trivial_minimizer} is optimal for $f = \chi_{R\,W_\phi}$.
		\end{remark}
		
		\begin{proposition}
			Let $f \in L^1(\R^n)$, $\lambda > 0$ and $u \in BV(\R^n)$. Then we have
			\[ E(u ; f, \lambda) = \int_{-\infty}^\infty E(\chi_{\Sigma_u(t)} ; \chi_{\Sigma_f(t)}, \lambda) ~\mathrm{d}t. \]
			
			Note that $\chi_{\Sigma_u(t)} \not \in L^1(\R^n)$ for $t < 0$, but the symmetric difference
			$\|\chi_{\Sigma_u(t)} - \chi_{\Sigma_f(t)} \|_{L^1(\R^n)} = | \Sigma_u(t) \symmdiff \Sigma_f(t) |$ is finite for almost every $t \in \R$.
			
			\begin{proof}
				We can follow the proof of \cite[Proposition 5.1]{Chan2005} since the anisotropic total variation satisfies the coarea formula.
			\end{proof}
		\end{proposition}		
		
		\begin{remark}
			If $f \in L^1(\R^n; \{ 0, 1\})$ is binary, it is a well known consequence of the coarea formula that there is a binary minimizer $u \in BV(\R^n; \{ 0, 1 \})$, see \cite[Theorem 5.2]{Chan2005}. In this special case the energy $E$ equals \[E(u) = \TV_\phi(u) + \lambda \| u -f \|_{L^2(\R^n)}^2,\] which is the anisotropic ROF model analyzed in \cite{Esedoglu2004}. Therefore everything that is proven in \cite[4.2]{Esedoglu2004} about the regularity of domains, that may appear as minimizer of the binary anisotropic ROF model, transfers to our situation.
			
		\end{remark}
		
		For $u \in BV(\R^n)$ we define the set (c.f. \cite{Choksi2011})
		\[
			\mathcal{V}(u) \coloneqq \{ v \in L^\infty(\R^n ; \R^n) \mid v \text{ satisfies (i) - (iii)} \} 
		\]
		with
		\begin{enumerate}[label=(\roman*)]
			\item $v \in -W_\phi$ almost everywhere in $\R^n$,
			\item has a weak divergence $\divergence v \in L^\infty(\R^n)$ and
			\item $\TV_\phi(u) = - \int_{\R^n} u \divergence v$.
		\end{enumerate}
		
		The following theorem is our main result.
		\begin{theorem}
			\label{thm:main}
			Let $f \in L^1(\R^n)$ and $\lambda > 0$. Then $u_0 \in BV(\R^n)$ is a minimizer of $E$ if and only if there is a $v \in \mathcal{V}(u_0)$ which additionally satisfies:
			\begin{enumerate}[label=(\alph*)]
				\item $\| \divergence v \|_{L^\infty(\R^n)} \le \lambda $,
				\item $\divergence v = \lambda$ almost everywhere in $\{ u_0 > f \}$ and
				\item $\divergence v = -\lambda$ almost everywhere in $\{ u_0 < f \}$.
			\end{enumerate}
			
			\begin{proof}
				 This proof is an extension of the proof of \cite[Lemma 3.1]{Choksi2011}.

				Let $v \in \mathcal{V}(u_0)$ be a function with the properties (a) -- (c). For every $h \in BV(\R^n)$ we have
				\begin{align*}
					E(u_0+h) &= \TV_\phi(u_0+h) + \lambda \| u_0 + h - f \|_{L^1(\R^n)} \\
					&\ge -\int_{\R^n} (u_0 + h) \divergence v + \lambda \| u_0 + h - f \|_{L^1(\R^n)} \\
					&= E(u_0) - \int_{\R^n} h \divergence v - \lambda \| u_0 - f \|_{L^1(\R^n)} + \lambda \|u_0 + h - f \|_{L^1(\R^n)} \\
					&= E(u_0) - \int_{\R^n} (h + u_0 - f) \divergence v + \lambda \|u_0 + h -f \|_{L^1(\R^n)} \\
					&\ge E(u_0).
				\end{align*}
				This proves that $u_0$ is a minimizer of $E$.
				
				Now let $u_0$ be a minimizer of $E$. We want to prove the existence of a $v \in \mathcal{V}(u_0)$ satisfying (a) -- (c).
				For that reason we define
				\[ \begin{array}{rl}
					F : BV(\R^n) \rightarrow \Mn ,~& u \mapsto Du, \\
					G : \Mn \rightarrow \R ,~& \mu \mapsto \int_{\R^n} \phi(\mathrm{d} \mu) \\
					H : BV(\R^n) \rightarrow \R ,~& u \mapsto \| u - f \|_{L^1(\R^n)}.
				\end{array} \]
				We have $E = G\circ F + \lambda H$.
				Since $F$, $G$ and $H$ are continuous we can apply subdifferential calculus (see \cite[I.5.6, I.5.7]{Ekeland1999}) and obtain
				\[ \partial E(u_0) = F^\ast \partial G(F u_0) + \lambda \partial H(u_0), \]
				where $F^\ast : (\Mn)^\ast\rightarrow BV(\R^n)^\ast$ is the transpose mapping of $F$.
				Since $u_0$ is a minimizer of $E$ we know that $0 \in \partial E(u_0)$. Therefore there are $\Phi \in \partial G(Fu_0) \subset (\Mn)^\ast$ and $\Psi \in \partial H(u_0) \subset BV(\R^n)^\ast$ with 
				\begin{equation}
					\label{eq:subdiff_zero}
					0 = F^\ast \Phi + \lambda \Psi.
				\end{equation}
				By definition we have
				\begin{equation}
					\label{eq:psi_norm}
					\langle \Psi, h \rangle \le H(h+u_0) - H(u_0) \le \| h \|_{L^1(\R^n)}
				\end{equation}
				for all $h \in BV(\R^n)$. Since $BV(\R^n)$ is dense in $L^1(\R^n)$ and $\Psi$ is a bounded operator on $BV(\R^n)$ with respect to the $L^1(\R^n)$ norm, we can extend $\Psi$ to $\Psi \in L^1(\R^n)^\ast = L^\infty(\R^n)$ with operator norm 
				\begin{equation}
					\label{eq:opnorm}
					\| \Psi \|_{L^\infty(\R^n ; \R^n)} \le 1.
				\end{equation}
				By continuity inequality \eqref{eq:psi_norm} remains valid for all $h \in L^1(\R^n)$. Testing this inequality with $\chi_{\{u_0 > f\}}(f-u_0)$ and $\chi_{\{u_0 < f\}}(f-u_0)$ we conclude
				\begin{equation}
					\label{eq:psi_exact}
					\begin{array}{ll}
						\Psi = 1 & \text{almost everywhere in } \{u_0 > f \} \text{ and } \\
						\Psi = -1 & \text{almost everywhere in } \{ u_0 < f \}.
					\end{array}
				\end{equation}				
				Plugging $Du_0$ and $-Du_0$ into the subdifferential inequality
				\[
					\langle \Phi, \mu - Du_0 \rangle \le G(\mu) - G(Du_0) ~\forall \mu \in \Mn
				\]
				gives
				\begin{equation}
					\label{eq:v_exact} \langle \Phi, Du_0 \rangle = \TV_\phi(u_0)
				\end{equation}
				and so
				\begin{equation}
					\label{eq:phi_norm}
					\langle \Phi, \mu \rangle \le G(\mu) ~\forall \mu \in \Mn.
				\end{equation}
				By restricting $\Phi$ to the space $L^1(\R^n ; \R^n)$, which is isometrically embedded in $\Mn$, we get $v = \restr{\Phi}{L^1(\R^n;\R^n)} \in (L^1(\R^n ; \R^n))^\ast = L^\infty(\R^n ; \R^n)$. Inequality \eqref{eq:phi_norm} yields
				\[
					\int_{\R^n} v \cdot g \le \int_{\R^n} \phi(g) ~ \forall g \in L^1(\R^n ; \R^n)
				\]
				and from \Cref{eq:hoelder_extremal} we can conclude that 
				\begin{equation}
					\label{eq:v_in_wphi}				
					v \in -W_\phi \text{ almost everywhere in } \R^n.
				\end{equation}
				Using equation \eqref{eq:subdiff_zero} we can deduce that the weak divergence of $v$ is given by
				\[ \divergence v = \lambda \Psi. \] 
				We conclude from \eqref{eq:v_exact} and \eqref{eq:v_in_wphi} that $v \in \mathcal{V}(u_0)$ and from \eqref{eq:opnorm} and \eqref{eq:psi_exact} the properties (a) -- (c).
			\end{proof}
		\end{theorem}
		
		\begin{remark}
			\Cref{thm:main} can also be obtained by the theory developed in \cite{Duval2011} if $\phi$ is even. We prefer to give a self-contained proof.
		\end{remark}
		
		\begin{corollary}\
			\begin{enumerate}
				\item We have that $u_0 \in BV(\R^n)$ is a minimizer of $E(\cdot ; u_0, \lambda)$ for $\lambda > 0$ if and only if there is a $v \in \mathcal{V}(u_0)$ satisfying $\| \divergence v \|_{L^\infty(\R^n)} \le \lambda$. If $\| \divergence v \|_{L^\infty(\R^n)} < \lambda$, then $u_0$ is the unique minimizer of $E(\cdot ; u_0, \lambda)$.
				\item If $u_0 \in BV(\R^n)$ is a minimizer of $E$ for some $f \in L^1(\R^n)$ and $\lambda >0$, then $u_0$ is a minimizer of $E(\cdot ; u_0, \lambda)$. In that sense the minimization of $E$ is an idempotent operation.
			\end{enumerate}
			
			\begin{proof}
			\begin{enumerate}
				\item The first part is a simple conclusion from \Cref{thm:main}. We can repeat the first part of the proof of \Cref{thm:main} to show that $u_0$ is the unique minimizer if $\| \divergence v \|_{L^\infty(\R^n)} < \lambda$.
				\item If $u_0$ is a minimizer of $E$ for some $f$, then we can apply \Cref{thm:main} to deduce the existence of a vector field $v \in \mathcal{V}(u_0)$ with $\| \divergence v \|_{L^\infty(\R^n)} \le \lambda$.
			\end{enumerate}
			\end{proof}
		\end{corollary}

		\begin{example}
			\label{ex}
			\begin{figure}[h!]
			\centering
			\begin{tikzpicture}[baseline=0, scale=3] 
				\draw[gray, dashed] (0,0) circle [radius=1cm];
				\draw[gray, dashed] (-0.9,-0.9) rectangle (0.9,0.9);

				\filldraw[darkgray] (0,0) circle(0.02cm);
 				\draw[darkgray] (0,0) -- node [right] {$h$} (0, 0.9)
 					-- 	node [below] {$s$} (0.43588989cm, 0.9cm)
 					-- node [right] {$1$} (0,0);
				
 				\draw[very thick] (0, 0.9cm) -- (0.43588989cm, 0.9cm) arc [radius=1, start angle=64.158, end angle=25.8419]
 					-- (0.9cm, -0.43588989cm) arc [radius=1, end angle=-64.158, start angle=-25.8419]
 					-- (-0.43588989cm, -0.9cm) arc [radius=1, start angle=244.158, end angle=205.8419]
 					-- (-0.9cm, (0.43588989cm) arc [radius=1, end angle=115.842, start angle=154.1581] -- cycle;
			\end{tikzpicture}
			\caption{The situation in \Cref{ex}.}
			\end{figure}
			In this example we consider the anisotropy $\phi = |\cdot|_1$ in dimension $n = 2$. The dual anisotropy is given by $\phi^\circ = |\cdot|_\infty$ and the Wulff shape is given by 
			\[ W_\phi = [-1,1]^2. \]
			We are going to calculate an optimal shape for a circle of radius $1$, which means $B = \{ |x|_2 \le 1 \}$ and $f = \chi_B$.
			It is known by \cite[Proposition 7.2.2]{Duval2011} that an optimal shape for an arbitrary bounded convex set can be expressed by an opening with a rescaled Wulff shape in the sense of morphology. In our case we get an optimal shape $U$, such that $u = \chi_U$ is a minimizer of $E$, either by $U = \emptyset$ or for $s = \frac{1}{\lambda}$ by
			\[
				U = B_s \coloneqq \left\{ x + sW_\phi \mid x \in B \text{ such that } x + sW_\phi \subset B \right\}.
			\]
			We can rewrite $B_s = B \cap \sqrt{1 - s^2}W_\phi = B \cap[-h,h]^2$ with $h = \sqrt{1 - s^2}$ if $h \in [\frac{1}{\sqrt{2}},1]$. We have
			\begin{align*}
				\TV_\phi (U) &= 8\sqrt{1 - \frac{1}{\lambda^2}}, \\
				|U| &= \pi - 4\left(\sin^{-1}\left(\frac{1}{\lambda}\right) - \frac{1}{\lambda}\sqrt{1 - \frac{1}{\lambda^2}}\right), \\
				E(\chi_U) &= 8h + 4 \lambda \left(\sin^{-1}\left(\frac{1}{\lambda}\right) - \frac{1}{\lambda}\sqrt{1 - \frac{1}{\lambda^2}}\right).
			\end{align*}
			From \cite[Proposition 7.2.2]{Duval2011} we conclude that $U$ is an optimal shape as long as 
			\[
				\frac{\TV_\phi(U)}{|U|} \ge \lambda,
			\]
			which is equivalent to
			\[
				4 \lambda \sin^{-1}\left(\frac{1}{\lambda}\right) + 4\sqrt{1 - \frac{1}{\lambda^2}} \ge \lambda \pi.
			\]
			This inequality is true as long as $\lambda \ge 2.4754\dotsc$.
			
			In the following we will prove that $U$ is an optimal shape for $\lambda > 2\sqrt{2}$. This is strictly weaker than the result we can deduce from \cite[Proposition 7.2.2]{Duval2011}, but the proof we give is shorter and hopefully more accessible. Furthermore this result is still strong enough to conclude that in the sense of \cite{Choksi2011} $E$ may be faithful to domains which are not clean 2D bar codes. For that reason we will construct a function $v$ which
			satisfies the conditions of \Cref{thm:main}.
			\begin{proof}
				Set
				\[
					w(x_1, x_2) \coloneqq \left\{ \begin{array}{cr}
						\min \{ 1, \max \{ -1, \frac{x_1}{s}\} \} & |x_2| \ge \frac{1}{\sqrt{2}} \\
						\min \{ 1, \max \{ -1, \sqrt{2}x_1 \} \} & |x_2| < \frac{1}{\sqrt{2}}
					\end{array}\right.
				\]
				and
				\[
					v(x) \coloneqq \left( \begin{array}{c}
					-w(x_1, x_2) \\
					-w(x_2, x_1)
					\end{array}\right).
				\]
				The construction implies $v(z) \cdot \nu(z) = -|\nu(z)|_1$ for $z \in \partial U$ and therefore
				\[
					\TV_{|\cdot|_1}(u) = \int_{\partial U} |\nu|_1 = -\int_{\R^n} u \divergence v. \]
				From \Cref{thm:main} we conclude, that $u$ is a minimizer of $E$.
			\end{proof}
		\end{example}
	
	\section*{Acknowledgments} This paper is an extension of parts of my master thesis, which i have finished in October 2016 at the Technische Universität Dortmund. I thank my advisor, Prof. Dr. Matthias Röger, for drawing my attention to this interesting topic and for his guidance during writing the thesis and this work.

	\bibliography{../literature}
\end{document}